\magnification=1200
\NoBlackBoxes
\loadeusm
\input BoxedEPS
\input amsppt.sty
\SetepsfEPSFSpecial
\topmatter
\title
Floating Body, Illumination Body, 
and Polytopal Approximation 
\endtitle
\author
Carsten Sch\"utt
\endauthor
\address
Mathematisches Seminar, 
Christian Albrechts Universit\"at,
D-24098 Kiel,
Germany
\endaddress
\address
Department of Mathematics,
Oklahoma State University,
Stillwater, Oklahoma 74078
\endaddress
\thanks
This paper was written while the author was visiting the MSRI
at Berkeley in the spring of 1996
\endthanks
\subjclass
52A22
\endsubjclass
\abstract
Let $K$ be a convex body in $\Bbb R^{d}$ and $K_{t}$ its
floating bodies. There is a polytope with at most
$n$ vertices that satisfies
$$
K_{t} \subset P_{n} \subset K
$$
where
$$
n \leq e^{16d} \frac{vol_{d}(K \setminus K_{t})}{t\ vol_{d}(B_{2}^{d})}
$$
Let $K^{t}$ be the illumination bodies of $K$ and $Q_{n}$ a
polytope that contains $K$ and has at most $n$ $d-1$-dimensional
faces. Then
$$
vol_{d}(K^{t} \setminus K) \leq cd^{4} vol_{d}(Q_{n} \setminus K)
$$
where 
$$
n \leq \frac{c}{dt} \ vol_{d}(K^{t} \setminus K)
$$
\endabstract
\endtopmatter
\newpage

\document

\heading
1. Introduction
\endheading

We investigate the approximation of a convex body $K$ in $\Bbb R^{d}$
by a polytope. We measure the approximation by the symmetric difference
metric. The symmetric difference metric between two convex bodies $K$
and $C$ is
$$
d_{S}(C,K)=vol_{d}((C \setminus K) \cup (K \setminus C))
$$
We study in particular two questions: How well can a convex body $K$ be 
approximated by a polytope $P_{n}$ that is contained in $K$ and has
at most $n$ vertices and how well can $K$ be approximated by a polytope
$Q_{n}$ that contains $K$ and has at most $n$ $d-1$-dimensional faces.
Macbeath [Mac] showed that the Euclidean Ball $B_{2}^{d}$ is an extremal
case: The approximation for any other convex body is better. We have for
the Euclidean ball
$$
c_{1}\ d\ vol_{d}(B_{2}^{d})n^{-\frac{2}{d-1}} \leq
d_{S}(P_{n},B_{2}^{d}) \leq c_{2}\ d\ vol_{d}(B_{2}^{d})n^{-\frac{2}{d-1}}
\tag 1.1
$$
provided that $n \geq (c_{3}\ d)^{\frac{d-1}{2}}$. The right hand inequality 
was first established by Bronshtein and Ivanov [BI] and Dudley [D$_{1}$,D$_{2}$].
Gordon, Meyer, and Reisner [GMR$_{1}$,GMR$_{2}$] gave a constructive proof 
for the same inequality. M\"uller
[M\"u] showed that random approximation gives the same estimate.
Gordon, Reisner, and Sch\"utt [GRS] established the left hand inequality.
Gruber [$\text{Gr}_{2}$] obtained an asymptotic formula. If a convex body 
$K$ in $\Bbb{R}^{d}$ has a $C^{2}$-boundary 
with everywhere positive curvature, then 
$$
\inf \{d_{S}(K,P_{n})|P_{n} \subset K\ \text{and $P_{n}$ \ has at most
n vertices} \}
$$
is asymptotically the same as
$$
\frac{1}{2} \text{del}_{d-1} \left(\int_{\partial K} \kappa(x)^{\frac{1}{d+1}}
d\mu(x)\right)^{\frac{d+1}{d-1}} 
(\frac{1}{n})^{\frac{2}{d-1}}
$$
where $del_{d-1}$ is a constant that is connected with Delone triangulations.
In this paper we are not concerned with asymptotic estimates, but with
uniform.
\par 
$Int(M)$ denotes the interior of a set $M$. $H(x,\xi)$ denotes the
hyperplane that contains $x$ and is orthogonal to $\xi$. $H^{+}(x,\xi)$
denotes the halfspace that contains the vector $x-\xi$, and
$H^{-}(x,\xi)$ the halfspace containing $x+\xi$.
$e_{i},i=1,\dots,d$ denotes the unit vector basis in $\Bbb R^{d}$.
$[A,B]$ is the convex hull of the sets $A$ and $B$.
The convex floating body $K_{t}$ of a convex body $K$ is the 
intersection of all halfspaces whose defining hyperplanes cut off
a set of volume $t$ from $K$.
\par
The illumination body $K^{t}$ of a convex body $K$ is [W] 
$$
\{x \in \Bbb R^{d}|\ vol_{d}([x,K] \setminus K) \leq t \ \}
$$
$K^{t}$ is a convex body. It is enough to show this for polytopes.
Let $F_{i}$ denote the faces of a polytope $P$, $\xi_{i}$ the outer normal
and $x_{i}$ an element of $F_{i}$. Then we have
$$
vol_{d}([x,P] \setminus P)=\frac{1}{d}\sum_{i=1}^{n}
\max\{0,<\xi_{i},x-x_{i}>\} vol_{d-1}(F_{i})
$$
The right hand side is a convex function.
\vskip 3mm

\heading
2. The Floating Body
\endheading

\proclaim{Theorem 2.1}
Let $K$ be a convex body in $\Bbb R^{d}$. Then we have for every
$t$, $0 \leq t \leq \frac{1}{4}e^{-4} vol_{d}(K)$, that there are $n \in \Bbb N$ with
$$
n \leq e^{16d}\ \frac{vol_{d}(K \setminus K_{t})}{t\ vol_{d}(B_{2}^{d})}
$$
and a polytope $P_{n}$ that has $n$ vertices and such that
$$
K_{t} \subset P_{n} \subset K
$$
\endproclaim
\vskip 3mm

We want to see what kind of asymptotic estimate we get for bodies with
smooth boundary from Theorem 1. We have [SW]
$$
vol_{d}(K \setminus K_{t}) \sim
t^{\frac{2}{d+1}} 
\frac{1}{2} \left(\frac{d+1}{vol_{d-1}(B_{2}^{d-1})} \right)^{\frac{2}{d+1}}
\int_{\partial K} \kappa(x)^{\frac{1}{d+1}}d\mu(x)
$$
$$
\sim t^{\frac{2}{d+1}} d \int_{\partial K} \kappa(x)^{\frac{1}{d+1}}d\mu(x)
$$
Since
$$
n \sim d^{\frac{d}{2}}\frac{1}{t} vol_{d}(K \setminus K_{t})
$$
we get
$$
vol_{d}(K \setminus K_{t}) \sim d \left(d^{\frac{d}{2}}\frac{1}{n}
vol_{d}(K \setminus K_{t}) \right)^{\frac{2}{d+1}}
\int_{\partial K} \kappa(x)^{\frac{1}{d+1}}d\mu(x)
$$
$$
vol_{d}(K \setminus K_{t})^{\frac{d-1}{d+1}} \sim d^{2}
n^{-\frac{2}{d+1}}\int_{\partial K} \kappa(x)^{\frac{1}{d+1}}d\mu(x)
$$
Thus we get
$$
vol_{d}(K \setminus P_{n}) \leq vol_{d}(K \setminus K_{t})
\sim d^{2}n^{-\frac{2}{d-1}}
\left( \int_{\partial K} \kappa(x)^{\frac{1}{d+1}}d\mu(x) \right)
^{\frac{d+1}{d-1}}
$$
In case that $K$ is the Euclidean ball we get
$$
vol_{d}(B_{2}^{d} \setminus P_{n}) \leq
cd^{2}n^{-\frac{2}{d-1}}vol_{d}(B_{2}^{d})
$$
where $c$ is an absolute constant.
If one compares this to the optimal result (1.1) one sees that there is
an additional factor $d$.
\par
The volume difference $vol_{d}(P)-vol_{d}(P_{t})$ for a polytope $P$ is of a 
much smaller order than for a convex body with smooth boundary. 
In fact, we have [S] that it is of the order $t|\ln t|^{d-1}$.
In [S] this has been used to get estimates for approximation of convex bodies
by polytopes.
\par
The same result as in Theorem 2.1 holds if we fix the number of (d-1)-dimensional 
faces instead of the number of vertices. This follows from the economic cap covering
for floating bodies [BL, Theorem 6]. The constants are not as good as in Theorem 2.1.
\vskip 3mm

The following lemmata are not new. They have usually been formulated for symmetric, 
convex bodies [B,H,MP].
\vskip 3mm

\proclaim{Lemma 2.2}
Let $K$ be a convex body in $\Bbb R^{d}$ and let $H(cg(K),\xi)$ be the
hyperplane passing through the center of gravity $cg(K)$ of $K$ and being
orthogonal to $\xi$.
Then we have for all $\xi \in \partial B_{2}^{d}$
\newline
(i)
$$
(1-\frac{1}{d+1})^{d} vol_{d}(K) \leq vol_{d}(K \cap H^{+}(cg(K),\xi))
\leq (1-(1-\frac{1}{d+1})^{d}) vol_{d}(K)
$$
(ii) for all hyperplanes $H$ in $\Bbb R^{d}$ that are parallel to 
$H(cg(K),\xi)$
$$
(1-\frac{1}{d+1})^{d-1} vol_{d-1}(K \cap H) \leq vol_{d-1}(K \cap H(cg(K),\xi))
$$
\endproclaim
\vskip 3mm

The sequence $(1-\frac{1}{d+1})^{d}$, $d=2,3,\dots$ is monotonely
decreasing. Indeed, by Bernoulli's inequality we have $1-\frac{1}{d}
\leq (1-\frac{1}{d^{2}})^{d}$, or $\frac{d-1}{d} \leq (\frac{d^{2}-1}{d^{2}})^{d}$. 
Therefore we get $(\frac{d}{d+1})^{d}
\leq(\frac{d-1}{d})^{d-1}$, which implies $(1-\frac{1}{d+1})^{d} \leq
(1-\frac{1}{d})^{d-1}$. 
\par
Therefore we get for the inequalities (i)
$$
\frac{1}{e} vol_{d}(K) \leq vol_{d}(K \cap H^{+}(cg(K),\xi))
\leq (1-\frac{1}{e}) vol_{d}(K)
\tag 2.1
$$
By the above $(1+\frac{1}{d})^{d}$ is a monotonely increasing sequence. Thus we
get $(1+\frac{1}{d})^{d-1} <e$. For (ii) we get
$$
vol_{d-1}(K \cap H) \leq e\ vol_{d-1}(K \cap H(cg(K),\xi))
\tag 2.2
$$

\demo{Proof}
(i) We can reduce the inequality to the case that $K$ is a cone
with a Euclidean ball of dimension $d-1$ as base. To see this we perform a 
Schwarz
symmetrization parallel to $H(cg(K),\xi)$ and denote the symmetrized body
by $S(K)$. The Schwarz symmetrization replaces a section parallel to 
$H(cg(K),\xi)$ by a $d-1$-dimensional Euclidean sphere of the same $d-1$-dimensional
volume. This does not change the volume of $K$ 
and $K \cap H^{+}(cg(K),\xi)$ and the center of gravity $cg(K)$ is still an element
of $H(cg(K),\xi)$. Now we consider the cone 
$$
[z,S(K) \cap H(cg(K),\xi)]
$$ 
such that
$$
vol_{d}([z,S(K) \cap H(cg(K),\xi)])=vol_{d}(K \cap H^{-}(cg(K),\xi))
$$
and such that $z$ is an element of the axis of symmetry of $S(K)$ and
of $H^{-}(cg(K),\xi)$. See figure 2.1.
$$
\tilde K=(K \cap H^{+}(cg(K),\xi)) \cup [z,S(K) \cap H(cg(K),\xi)]
$$ 
is a convex set such that $vol_{d}(K)=vol_{d}(\tilde K)$ and such that
the center of gravity $cg(\tilde K)$ of $\tilde K$ is contained in
$[z,S(K) \cap H(cg(K),\xi)]$. Thus
$$
vol_{d}(\tilde K \cap  H^{+}(cg(\tilde K),\xi)) \geq 
vol_{d}(\tilde K \cap  H^{+}(cg(K),\xi)) =vol_{d}(K \cap H^{+}(cg(K),\xi))
$$
We apply a similar argument to the set $S(K) \cap H^{+}(cg(K),\xi)$ and show
that we may assume that $S(K)$ is a cone with $z$ as its vertex. 
Thus we may assume that
$$
K=[(0,\dots,0,1),\{(x_{1},\dots,x_{d-1},0)| \sum_{i=1}^{d-1}|x_{i}|^{2} \leq 1\}]
\quad  \text{and} \quad
\xi=(0,\dots,0,1)
$$
\vskip 3mm

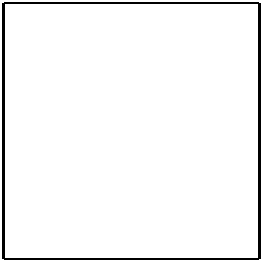
\BoxedEPSF{convex2 scaled 800}
\vskip 5mm

Then
$$
vol_{d}(K)=\frac{1}{d}vol_{d-1}(B_{2}^{d-1})
$$
and
$$
\frac{1}{vol_{d}(K)}\int_{K} x_{d}dx_{d}=
d\int_{0}^{1}t(1-t)^{d-1}dt=d\int_{0}^{1}(1-s)s^{d-1}ds=\frac{1}{d+1}
$$
We obtain that
$$
vol_{d}(K \cap H^{-}(cg(K),(0,\dots,0,1))=(1-\frac{1}{d+1})^{d}\ vol_{d}(K)
$$
(ii) Let $H$ be a hyperplane that is parallel to $H(cg(K),\xi)$ and such that
$vol_{d-1}(K \cap H)>vol_{d-1}(K \cap H(cg(K),\xi))$. Otherwise there is
nothing to prove. We apply a Schwarz 
symmetrization parallel to $H(cg(K),\xi)$ to $K$. The symmetrized body is
denoted by $S(K)$. Let $z$ be the element 
of the axis of symmetry of $S(K)$ such that
$$
[z,S(K) \cap H] \cap H(cg(K),\xi)=S(K) \cap H(cg(K),\xi)
$$
Since $vol_{d-1}(K \cap H)>vol_{d-1}(K \cap H(cg(K),\xi))$ there is such a
$z$. 
We may assume that $H^{+}(cg(K),\xi)$ is the half space containing $z$. 
Then we have
$$
\aligned
[z,S(K) \cap H] \cap H^{-}(cg(K),\xi) & \subset S(K) \cap H^{-}(cg(K),\xi) \\
[z,S(K) \cap H] \cap H^{+}(cg(K),\xi) & \supset S(K) \cap H^{+}(cg(K),\xi)
\endaligned
$$
Therefore we have that
$$
cg([z,S(K) \cap H]) \in H^{+}(cg(K),\xi)
$$
Therefore, if $h_{cg}$ denotes the distance of $z$ to $H(cg(K),\xi)$
and $h$ the distance of $z$ to $H$, we get as in the proof of (i) that
$$
h_{cg} \geq h(1-\frac{1}{d+1})
$$
Thus we get
$$
vol_{d-1}(K \cap H(cg(K),\xi))=
vol_{d-1}(S(K) \cap H(cg(K),\xi)) \geq
$$
$$
(1-\frac{1}{d+1})^{d-1}vol_{d-1}(S(K) \cap H) =
(1-\frac{1}{d+1})^{d-1}vol_{d-1}(K \cap H)
$$
\enddemo
\qed
\vskip 3mm

\proclaim{Lemma 2.3}
Let $K$ be a convex body in $\Bbb R^{d}$ and let $\Theta(\xi)$
be the infimum of
all numbers $t$, $0<t$, such that
$$
vol_{d-1}(K \cap H(cg(K),\xi)) \geq
e\ vol_{d-1}(K \cap H(cg(K)+t\xi,\xi))
$$
Then we have
$$
\frac{1}{2e^{3}}vol_{d}(K) \leq \Theta(\xi)vol_{d-1}(K \cap H(cg(K),\xi))
\leq e\ vol_{d}(K)
$$ 
\endproclaim
\vskip 3mm

\demo{Proof}
The right hand inequality follows from Fubini\rq s theorem and Brunn-Minkowski\rq s
theorem. Now we verify the left hand inequality. We consider first
the case that we have for $t$, $t> \Theta(\xi)$, 
$$
K \cap H(cg(K)+t\xi,\xi) = \emptyset
$$
Then we have by (2.1) and (2.2)
$$
\aligned
\frac{1}{e}vol_{d}(K) \leq & vol_{d}(K \cap H^{+}(cg(K),\xi))  \\
= & \int_{0}^{\Theta(\xi)}vol_{d-1}(K \cap H(cg(K)+t\xi,\xi))dt
\leq e\ \Theta(\xi)\ vol_{d-1}(H(cg(K),\xi))
\endaligned
$$
If for some $t$, $t> \Theta(\xi)$, we have $ K \cap H(cg(K)+t\xi,\xi) \ne \emptyset$
then we have
$$
vol_{d-1}(K \cap H(cg(K),\xi)) =
e\ vol_{d-1}(K \cap H(cg(K)+\Theta(\xi)\xi,\xi))
$$ 
We perform a Schwarz 
symmetrization parallel to $H(cg(K),\xi)$. We consider the cone
$$
[z,S(K) \cap H(cg(K),\xi)]
$$ 
such that $z$ is an element of the axis of symmetry of $S(K)$ and such that
$$
[z,S(K) \cap H(cg(K),\xi)] \cap H(cg(K)+\Theta(\xi)\xi,\xi)=
S(K) \cap H(cg(K)+\Theta(\xi)\xi,\xi)
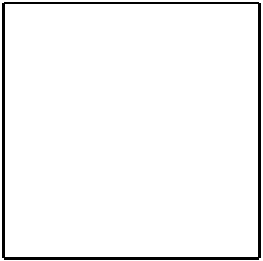
$$
\vskip 3mm
\BoxedEPSF{convex3 scaled 800}
\vskip 5mm

Let $ H^{+}(cg(K),\xi)$ and $ H^{+}(cg(K)+\Theta(\xi)\xi,\xi)$ be the
half spaces that contain $z$. Then we get by convexity
$$
\aligned
[z,S(K) \cap H(cg(K),\xi)] \cap & H^{+}(cg(K)+\Theta(\xi)\xi,\xi)  \\
& \supset S(K) \cap H^{+}(cg(K)+\Theta(\xi)\xi,\xi)
\endaligned
\tag 2.3
$$
We get by (2.1)
$$
\frac{1}{e}vol_{d}(K) \leq vol_{d}(K \cap H^{+}(cg(K),\xi)) = 
$$
$$
\aligned
vol_{d}(K \cap H^{+}(cg(K),\xi) \cap H^{-}(cg(K)+\Theta(\xi)\xi,\xi))+  \\
vol_{d}(K \cap H^{+}(cg(K)+ & \Theta(\xi)\xi,\xi)) =
\endaligned   
$$
$$
\aligned
vol_{d}(S(K) \cap H^{+}(cg(K),\xi) \cap H^{-}(cg(K)+\Theta(\xi)\xi,\xi))+  \\
vol_{d}(S(K) \cap H^{+}(cg(K)+ & \Theta(\xi)\xi,\xi))
\endaligned    
$$
By the hypothesis of the lemma we have for all $s$ with $0 \leq s \leq \Theta(\xi)$
$$
vol_{d-1}(K \cap H(cg(K),\xi)) \leq e\ vol_{d-1}(K \cap H(cg(K)+s\xi,\xi))
$$
Using this and (2.2) we estimate the first summand. The second summand is
estimated by using (2.3). Thus the above expression is not greater than
$$
\aligned
e^{2}\ vol_{d}([z,S(K) \cap H(cg(K),\xi)] 
\cap H^{-}(cg(K)+\Theta(\xi)\xi,\xi))+    \\
vol_{d}([z,S(K) \cap H(cg(K),\xi)] \cap H^{+} & (cg(K)+\Theta(\xi)\xi,\xi))
\endaligned 
$$
By an elementary computation for the volume of a cone we get that the latter
expression is smaller than
$$
2e^{2}
vol_{d}([z,S(K) \cap H(cg(K),\xi)] \cap H^{-}(cg(K)+\Theta(\xi)\xi,\xi))
$$
We use (2.2) again and get that the above expression is smaller than
$$
2e^{3}\Theta(\xi) vol_{d-1}(K \cap H(cg(K),\xi))
$$
\enddemo
\qed
\vskip 3mm

\proclaim{Lemma 2.4}
Let $K$ be a convex body in $\Bbb R^{d}$. Then there is a linear
transform $T$ with $\det(T)=1$ so that we have for all 
$\xi \in \partial B_{2}^{d}$
$$
\int_{T(K)}|<x,\xi>|^{2}dx=\frac{1}{d}
\int_{T(K)} \sum_{i=1}^{d} |<x,e_{i}>|^{2}dx
$$
\endproclaim
\vskip 3mm

We say that a convex body is in an isotropic position if the linear
transform $T$ in Lemma 2.4 can be chosen to be the identity. See [B,H].

\demo{Proof}
We claim that there is a orthogonal transform $U$ such that we have
for all $i,j=1, \dots,d$ with $i \ne j$,
$$
\int_{U(K)} <x,e_{i}><x,e_{j}>dx=0
$$
Clearly, the matrix
$$
(\int_{K} <x,e_{i}><x,e_{j}>dx)_{i,j=1}^{d}
$$
is symmetric. Therefore there is an orthogonal $d \times d$-matrix $U$
so that
$$
U(\int_{K} <x,e_{i}><x,e_{j}>dx)_{i,j=1}^{d}U^{t}
$$
is a diagonal matrix. We have
$$
U(\int_{K} <x,e_{i}><x,e_{j}>dx)_{i,j=1}^{d}U^{t}
=(\int_{K} \sum_{i,j=1}^{d} u_{l,i}<x,e_{i}><x,e_{j}>u_{k,j} dx)_{l,k=1}^{d}
$$
$$
=(\int_{K}<x,U^{t}(e_{l})><x,U^{t}(e_{k})>dx)_{l,k=1}^{d}
=(\int_{U(K)}<y,e_{l}><y,e_{k}>dy)_{l,k=1}^{d}
$$
So the latter matrix is a diagonal matrix. All the diagonal elements
are strictly positive. This argument is repeated with a diagonal matrix so
that the diagonal elements turn out to be equal. Therefore there is a matrix
$T$ with $\det T=1$ such that
$$
\int_{T(K)} <x,e_{i}><x,e_{j}>dx= \left\{
\aligned
0  & \hskip 42mm \text{if} \quad i \ne j \\
\frac{1}{d} & \int_{T(K)} \sum_{j=1}^{d}|<x,e_{j}>|^{2}dx
\quad  \text{if} \quad i=j
\endaligned
\right .
$$
From this the lemma follows.
\enddemo
\qed
\vskip 3mm

\proclaim{Lemma 2.5} Let K be a convex body in $\Bbb R^{d}$ that is in an
isotropic position and whose center of gravity is at the origin. 
Then we have for all $\xi \in \partial B_{2}^{d}$
$$
\frac{1}{24e^{10}} vol_{d}(K)^{3} \leq
vol_{d-1}(K \cap H(cg(K),\xi))^{2}\frac{1}{d}
\int_{K} \sum_{i=1}^{d} |<x,e_{i}>|^{2}dx
\leq 6\ e^{3}\ vol_{d}(K)^{3}
$$
\endproclaim
\vskip 3mm

\demo{Proof}
By Lemma 2.4 we have for all $\xi \in \partial B_{2}^{d}$
$$
\frac{1}{d}
\int_{K} \sum_{i=1}^{d} |<x,e_{i}>|^{2}dx=
\int_{K}|<x,\xi>|^{2}dx
$$ 
By Fubini`s theorem we get that this equals
$$
\int_{-\infty}^{\infty}t^{2}\ 
vol_{d-1}(K \cap H(t\xi,\xi))dt \geq
\int_{0}^{\Theta(\xi)}t^{2}\ 
vol_{d-1}(K \cap H(t\xi,\xi))dt 
$$
where $\Theta(\xi)$ is as defined in Lemma 2.3.
By the definition of $\Theta(\xi)$ the above expression is greater than
$$
\frac{1}{e}vol_{d-1}(K \cap H(cg(K),\xi))\int_{0}^{\Theta(\xi)}t^{2}\ 
dt \geq
\frac{1}{3e}\Theta(\xi)^{3} vol_{d-1}(K \cap H(cg(K),\xi)) 
$$
By Lemma 2.3 this is greater than
$$
\frac{1}{24e^{10}}\frac{vol_{d}(K)^{3}}
{vol_{d-1}(K \cap H(cg(K),\xi))^{2}}
$$
Now we show the right hand inequality. By Lemma 2.4 we have
$$
\frac{1}{d}\int_{K} \sum_{i=1}^{d} |<x,e_{i}>|^{2}dx=
\int_{K}|<x,\xi>|^{2}dx=
\int_{-\infty}^{\infty}t^{2}\ 
vol_{d-1}(K \cap H(t\xi,\xi))dt=
$$
$$
\int_{0}^{\Theta(\xi)}t^{2}\ 
vol_{d-1}(K \cap H(t\xi,\xi))dt
+\int_{\Theta(\xi)}^{\infty}t^{2}\ 
vol_{d-1}(K \cap H(t\xi,\xi))dt+ 
$$
$$
\int^{0}_{\Theta(-\xi)}t^{2}\ 
vol_{d-1}(K \cap H(t\xi,\xi))dt
+\int^{\Theta(-\xi)}_{-\infty}t^{2}\ 
vol_{d-1}(K \cap H(t\xi,\xi))dt
$$
By (2.2) this is not greater than
$$
\frac{e}{3}\Theta(\xi)^{3}vol_{d-1}(K \cap H(cg(K),\xi))+
\int_{\Theta(\xi)}^{\infty}t^{2}\ 
vol_{d-1}(K \cap H(t\xi,\xi))dt+ 
$$
$$
\frac{e}{3}\Theta(-\xi)^{3}vol_{d-1}(K \cap H(cg(K),\xi))+
\int^{\Theta(-\xi)}_{-\infty}t^{2}\ 
vol_{d-1}(K \cap H(t\xi,\xi))dt 
$$
The integrals can be estimated by
$$
2\ \Theta(\xi)^{3}vol_{d-1}(K \cap H(cg(K),\xi))
\quad \text{and}\quad
2\ \Theta(-\xi)^{3}vol_{d-1}(K \cap H(cg(K),\xi))
$$
respectively. We treat here only the case $\xi$, the case $-\xi$ is treated in the
same way. If the integral equals $0$ then there is nothing to show. If the integral
does not equal $0$ then we have
$$
vol_{d-1}(K \cap H(cg(K),\xi)) =
e\ vol_{d-1}(K \cap H(cg(K)+\Theta(\xi)\xi,\xi))
$$ 
We consider the Schwarz symmetrization $S(K)$ of $K$ 
with respect to the plane $H(cg(K),\xi)$. We consider the cone $C$ that
is generated by the Euclidean spheres $S(K) \cap H(cg(K),\xi)$ and
$S(K) \cap H(cg(K)+\Theta(\xi)\xi,\xi)$. We have that
$$
S(K) \cap H^{+}(cg(K)+\Theta(\xi)\xi,\xi) \subset C
$$
and that the height of $C$ is equals
$$
\frac{\Theta(\xi)}{1-e^{-\frac{1}{d-1}}}
$$
Since $(1+\frac{1}{d-1})^{d-1}<e$ we have $1-e^{-\frac{1}{d-1}}>\frac{1}{d}$.
Thus the height of the cone $C$ is less than $d\ \Theta(\xi)$. Thus we get
for all $t$ with $\Theta(\xi) \leq t \leq d\ \Theta(\xi)$
$$
vol_{d-1}(K \cap H(cg(K)+t\xi,\xi) \leq (1-\frac{t}{d\Theta(\xi)})^{d-1}
vol_{d-1}(K \cap H(cg(K),\xi))
$$
Now we get
$$
\align
\int_{\Theta(\xi)}^{\infty}t^{2}\ 
vol_{d-1} & (K \cap H(t\xi,\xi))dt \leq           \\
& \int_{\Theta(\xi)}^{d\ \Theta(\xi)}t^{2}\ (1-\frac{t}{d\Theta(\xi)})^{d-1}
vol_{d-1}(K \cap H(cg(K),\xi)) dt \leq                      \\ 
& vol_{d-1}(K \cap H(cg(K),\xi)) (d\ \Theta(\xi))^{3} 
\int_{0}^{1}s^{2}(1-s)^{d-1}ds=                             \\
& vol_{d-1}(K \cap H(cg(K),\xi))(d\ \Theta(\xi))^{3} 
\frac{2}{d(d+1)(d+2)} \leq                                  \\
& 2\ vol_{d-1}(K \cap H(cg(K),\xi))\Theta(\xi)^{3}
\endalign
$$
Therefore we get
$$
\frac{1}{d}\int_{K} \sum_{i=1}^{d} |<x,e_{i}>|^{2}dx \leq
(\frac{e}{3}+2) (\Theta(\xi)^{3}+\Theta(-\xi)^{3})vol_{d-1}(K \cap H(cg(K),\xi))
$$
Now we apply Lemma 2.3 and get
$$
2(\frac{e}{3}+2) e^{3}\frac{vol_{d}(K)^{3}}{vol_{d-1}(K \cap H(cg(K),\xi))^{2}}
$$
\enddemo
\qed
\vskip 3mm

\proclaim{Lemma 2.6}
Let $K$ be a convex body in $\Bbb R^{d}$ such that the origin is an 
element of $K$. Then we have
$$
\frac{1}{d}
\int_{K} \sum_{i=1}^{d} |<x,e_{i}>|^{2}dx \geq
\frac{d^{\frac{2}{d}}}{d+2} vol_{d-1}(\partial B_{2}^{d})^{-\frac{2}{d}}
vol_{d}(K)^{\frac{d+2}{d}}
$$
\endproclaim
\vskip 3mm

\demo{Proof}
Let $r(\xi)$ be the distance of the origin to the boundary of $K$
in direction $\xi$. By passing to spherical coordinates we get
$$
\frac{1}{d}\int_{K} \sum_{i=1}^{d} |<x,e_{i}>|^{2}dx =
\frac{1}{d}\int_{\partial B_{2}^{d}}\int_{0}^{r(\xi)} \rho^{d+1} d\rho d\xi=
\frac{1}{d(d+2)}\int_{\partial B_{2}^{d}}r(\xi)^{d+2}d\xi
$$
By H\"older`s inequality we get that the above expression is greater than
$$
\frac{vol_{d-1}(\partial B_{2}^{d})}{d(d+2)}
\left(\frac{1}{vol_{d-1}(\partial B_{2}^{d})}
\int_{\partial B_{2}^{d}}r(\xi)^{d}d\xi\right)^{\frac{d+2}{d}}=
\frac{d^{\frac{2}{d}}}{d+2} vol_{d-1}(\partial B_{2}^{d})^{-\frac{2}{d}}
vol_{d}(K)^{\frac{d+2}{d}}
$$
\enddemo
\qed
\vskip 3mm

The following lemma can be found in [MP]. It is formulated there for
the case of symmetric convex bodies.

\proclaim{Lemma 2.7}
Let $K$ be a convex body in $\Bbb R^{d}$ such that the origin coincides with
the center of gravity of $K$ and such that $K$ 
is in an isotropic position. 
Then we have 
$$
B_{2}^{d}(cg(K),\frac{1}{24e^{5}\sqrt{\pi}}vol_{d}(K)^{\frac{1}{d}}) 
\subset K_{\frac{1}{4e^{4}}vol_{d}(K)}
$$
\endproclaim
\vskip 3mm

\demo{Proof}
As in Lemma 2.3 let $\Theta(\xi)$ be the infimum of
all numbers $t$ such that
$$
vol_{d-1}(K \cap H(cg(K),\xi)) \geq
e\ vol_{d-1}(K \cap H(cg(K)+t\xi,\xi))
$$ 
By Lemma 2.3 we have
$$
\Theta(\xi) \geq \frac{1}{2e^{3}}
\frac{vol_{d}(K)}{vol_{d-1}(K \cap H(cg(K),\xi))}
$$
By Lemma 2.5 we get 
$$
\Theta(\xi) \geq
\frac{1}{2e^{3}\sqrt{6}e^{\frac{3}{2}}}\left( \frac{1}{vol_{d}(K)}
\frac{1}{d}\int_{K}\sum_{i=1}^{d}|<x,e_{i}>|^{2}dx \right)^{\frac{1}{2}}
$$
We have 
$$
vol_{d}(B_{2}^{d})= \frac{\pi^{\frac{d}{2}}}{\Gamma(\frac{d}{2}+1)} 
\leq 
\frac{\pi^{\frac{d-1}{2}}(2e)^\frac{d}{2}}{d^{\frac{d+1}{2}}}
$$
and thus
$$
vol_{d}(B_{2}^{d})^{\frac{1}{d}} \leq
\sqrt{\frac{2\pi e}{d}}
$$
Therefore we get by Lemma 2.6 
$$
\Theta(\xi) \geq
\frac{1}{2e^{3}\sqrt{6}e^{\frac{3}{2}}}
\frac{d^{\frac{1}{d}}}{\sqrt{d+2}}
\left(\frac{vol_{d}(K)}{vol_{d-1}(\partial B_{2}^{d})}\right)^{\frac{1}{d}}
\geq \frac{1}{12e^{5}\sqrt{\pi}}\ vol_{d}(K)^{\frac{1}{d}} 
$$
On the other hand, we have
$$
vol_{d}(K \cap H^{-}(cg(K)+\frac{\Theta(\xi)}{2}\xi,\xi)) \geq
\int_{\frac{\Theta(\xi)}{2}}^{\Theta(\xi)}
vol_{d-1}(K \cap H(cg(K)+t\xi,\xi))dt
$$
where $H^{-}(cg(K)+\frac{\Theta(\xi)}{2}\xi,\xi)$ is the half space not
containing the origin. 
By the definition of $\Theta(\xi)$ this expression is greater than
$$
\frac{\Theta(\xi)}{2e}vol_{d-1}(K \cap H(cg(K),\xi))
$$
By Lemma 2.3 we get that this is greater than
$$
\frac{1}{4e^{4}}vol_{d}(K)
$$
Therefore, every hyperplane that has distance
$$
\frac{1}{24e^{5}\sqrt{\pi}} vol_{d}(K)^{\frac{1}{d}}
$$
from the center of gravity cuts off a set of volume greater than
$\frac{1}{4e^{4}}vol_{d}(K)$. 
\enddemo
\qed
\vskip 5mm

\demo{Proof of Theorem 2.1}
We are choosing the vertices $x_{1}, \dots, x_{n} \in \partial K$
of the polytope $P_{n}$. $N(x_{k})$ denotes the normal to $\partial K$
at $x_{k}$. $x_{1}$ is chosen arbitrarily. Having chosen
$x_{1}, \dots, x_{k-1}$ we choose $x_{k}$ such that
$$
\{x_{1}, \dots, x_{k-1} \} \cap Int(K \cap H^{-}(x_{k}-\Delta_{k}N(x_{k}),
N(x_{k}))= \emptyset
$$
where $\Delta_{k}$ is determined by
$$
vol_{d}(K \cap H^{-}(x_{k}-\Delta_{k}N(x_{k}),N(x_{k})))=t
$$
It could be that the hyperplane $H(x_{k}-\Delta_{k}N(x_{k}),N(x_{k}))$ is not
tangential to the floating body $K_{t}$, but this does not affect the
computation.
We claim that this process terminates for some $n$ with
$$
n \leq e^{16d}\frac{vol_{d}(K \setminus K_{t})}{t\ vol_{d}(B_{2}^{d})}
\tag 2.4
$$
This claim proves the theorem: If we cannot choose another $x_{n+1}$,
then there is no cap of volume $t$ that does not contain an element of the 
polytope $P_{n}=[x_{1}, \dots, x_{n}]$. By the theorem of Hahn-Banach
we get $K_{t} \subset P_{n}$. We show now the claim. We put
$$
\aligned
S_{n} &=K \cap H^{-}(x_{n}-\Delta_{n}N(x_{n}),N(x_{n}))  \\
S_{k} &=K \cap \left(\bigcap_{i=k+1}^{n}H^{+}(x_{i}-\Delta_{i}N(x_{i}),
N(x_{i})) \right) \cap H^{-}(x_{k}-\Delta_{k}N(x_{k}),N(x_{k}))
\endaligned
\tag 2.5
$$
for $k=1, \dots,n-1$.
We have for $k \ne l$ that
$$
vol_{d}(S_{k} \cap S_{l})=0
$$
Let $k <l <n$. Then we have
$$
\aligned
S_{k} \cap S_{l} & = 
K \cap \left(\bigcap_{i=k+1}^{n}H^{+}(x_{i}-\Delta_{i}N(x_{i}),
N(x_{i}))\right) \cap H^{-}(x_{k}-\Delta_{k}N(x_{k}),N(x_{k}))      \\
& \cap
K \cap \left(\bigcap_{i=l+1}^{n}H^{+}(x_{i}-\Delta_{i}N(x_{i}),
N(x_{i}))\right) \cap H^{-}(x_{l}-\Delta_{l}N(x_{l}),N(x_{l}))        \\
& \subset H^{+}(x_{l}-\Delta_{l}N(x_{l}),N(x_{l}))
\cap H^{-}(x_{l}-\Delta_{l}N(x_{l}),N(x_{l}))      \\
& =H(x_{l}-\Delta_{l}N(x_{l}),N(x_{l}))
\endaligned
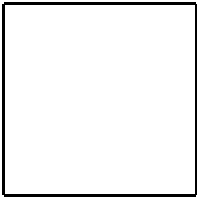
$$
\vskip 3mm
\BoxedEPSF{convex4 scaled 600}
\vskip 5mm

Thus we have
$$
vol_{d}(S_{k} \cap S_{l}) \leq
vol_{d}(H(x_{l}-\Delta_{l}N(x_{l}),N(x_{l})))=0
\tag 2.6
$$
The case $k<l=n$ is shown in the same way. 
We have for $k=1,\dots,n-1$
$$
\aligned
S_{k} & =K \cap \left(\bigcap_{i=k+1}^{n}H^{+}(x_{i}-\Delta_{i}N(x_{i}),
N(x_{i}))\right) \cap H^{-}(x_{k}-\Delta_{k}N(x_{k}),N(x_{k}))   \\
& \supset [x_{k},K_{t}] \cap H^{-}(x_{k}-\Delta_{k}N(x_{k}),N(x_{k}))  \\
& \supset [x_{k},(K \cap H^{-}(x_{k}-\tilde \Delta_{k}N(x_{k}),N(x_{k}))_{t}]
\cap H^{-}(x_{k}-\Delta_{k}N(x_{k}),N(x_{k}))
\endaligned
$$
where $\tilde \Delta_{k}$ is determined by
$$
vol_{d}(K \cap H^{-}(x_{k}-\tilde \Delta_{k}N(x_{k}),N(x_{k})))=4e^{4}t
$$
By Lemma 2.7 there is an ellipsoid $\eusm E$ contained in
$(K \cap H^{-}(x_{k}-\tilde \Delta_{k}N(x_{k}),N(x_{k})))_{t}$ whose center is
$cg(K \cap H^{-}(x_{k}-\tilde \Delta_{k}N(x_{k}),N(x_{k})))$
and that has volume
$$
vol_{d}(\eusm E)=\frac{4e^{4}}{(24e^{5}\sqrt{\pi})^{d}}\ t\ vol_{d}(B_{2}^{d})
$$
Since $(K \cap H^{-}(x_{k}-\tilde \Delta_{k}N(x_{k}),N(x_{k})))_{t}$ is contained
in $K_{t}$, $\eusm E$ is contained in $K_{t}$.
Thus
$$
S_{k} \supset [x_{k},\eusm E] \cap H^{-}(x_{k}- \Delta_{k}N(x_{k}),N(x_{k}))
$$
We claim now that $[x_{k},\eusm E] \cap H^{-}(x_{k}- \Delta_{k}N(x_{k}),N(x_{k}))$
contains an ellipsoid $\tilde \eusm E$ such that
$$
vol_{d}(\tilde \eusm E)=\frac{4e^{4}}{(24e^{5}\sqrt{\pi})^{d}}
\frac{1}{(4e^{5})^{d}}\ t\ vol_{d}(B_{2}^{d})
$$
and consequently
$$
vol_{d}(S_{k}) \geq
\frac{4e^{4}}{(24e^{5}\sqrt{\pi})^{d}}
\frac{1}{(4e^{5})^{d}}\ t\ vol_{d}(B_{2}^{d})=
\frac{4e^{4}}{(96e^{10}\sqrt{\pi})^{d}}\ t\ vol_{d}(B_{2}^{d})
\tag 2.7
$$ 
For this we have to see that
$\tilde \Delta_{k} \leq 4e^{5}\ \Delta_{k}$. By the assumption 
$t \leq \frac{1}{4}e^{-5}vol_{d}(K)$ we get that
$$
vol_{d}(K \cap H^{-}(x_{k}-\tilde \Delta_{k}N(x_{k}),N(x_{k}))) \leq
\frac{1}{e}\ vol_{d}(K)
$$
Therefore we get by (2.1) that
$cg(K) \in H^{+}(x_{k}-\tilde \Delta_{k}N(x_{k}),N(x_{k}))$. We consider 
two cases. If
$$
vol_{d-1}(K \cap H(x_{k}-\tilde \Delta_{k}N(x_{k}),N(x_{k}))) \leq
vol_{d-1}(K \cap H(x_{k}- \Delta_{k}N(x_{k}),N(x_{k})))
$$
then we have for all $t$, $\Delta_{k} \leq t \leq \tilde \Delta_{k}$, 
by the theorem of Brunn-Minkowski 
$$
\aligned
vol_{d-1}(K \cap H(cg(K),N(x_{k})) & \leq
vol_{d-1}(K \cap H(x_{k}-\tilde \Delta_{k}N(x_{k}),N(x_{k})))    \\  
& \leq vol_{d-1}(K \cap H(x_{k}- tN(x_{k}),N(x_{k})))            \\ 
\endaligned
\tag 2.8
$$
We get by (2.2)
$$
\Delta_{k} \geq \frac{t}{e\ vol_{d-1}(K \cap H(cg(K),N(x_{k})))}
$$
By (2.8)
$$
\align
(\tilde\Delta_{k} & -\Delta_{k}) vol_{d-1}(K \cap H(cg(K),N(x_{k}))) \leq   \\
& vol_{d}(K \cap H^{-}(x_{k}-\tilde \Delta_{k}N(x_{k}),N(x_{k})))-
vol_{d}(K \cap H^{-}(x_{k}-\Delta_{k}N(x_{k}),N(x_{k})))
\endalign
$$
This implies
$$
\tilde\Delta_{k}-\Delta_{k} \leq
\frac{(4e^{4}-1)t}{vol_{d-1}(K \cap H(cg(K),N(x_{k})))}
$$
Therefore we get
$$
\tilde\Delta_{k} \leq \frac{(4e^{4}-1)t}{vol_{d-1}(K \cap H(cg(K),N(x_{k})))}
+ \Delta_{k}
\leq 4e^{5}\ \Delta_{k}
$$
If
$$
vol_{d-1}(K \cap H(x_{k}- \Delta_{k}N(x_{k}),N(x_{k}))) \leq
vol_{d-1}(K \cap H(x_{k}-\tilde \Delta_{k}N(x_{k}),N(x_{k})))
$$
then by the theorem of Brunn-Minkowski we have for all $t$, 
$0 \leq t \leq \Delta_{k}$, and all $s$, $\Delta_{k} \leq s
\leq \tilde \Delta_{k}$,
$$
\aligned
vol_{d-1}(K \cap H(x_{k}-tN(x_{k}),N(x_{k}))) & 
\leq vol_{d-1}(K \cap H(x_{k}-\Delta_{k} N(x_{k}),N(x_{k})))   \\ 
& \leq vol_{d-1}(K \cap H(x_{k}-sN(x_{k}),N(x_{k})))
\endaligned
$$
We get
$$
\Delta_{k} \geq \frac{t}{vol_{d-1}(K \cap H(x_{k}-\Delta_{k} N(x_{k}),N(x_{k})))}
$$
and
$$
\tilde \Delta_{k}-\Delta_{k} \leq
\frac{(4e^{4}-1)t}{vol_{d-1}(K \cap H(x_{k}-\Delta_{k} N(x_{k}),N(x_{k})))}
$$
Therefore we get
$$
\tilde \Delta_{k} \leq
\frac{(4e^{4}-1)t}{vol_{d-1}(K \cap H(x_{k}-\Delta_{k} N(x_{k}),N(x_{k}))}+
\Delta_{k} \leq 4e^{4} \Delta_{k}
$$
We have verified (2.7). From (2.6) and (2.7) we get
$$
vol_{d}(K \setminus K_{t}) \geq vol_{d}(\bigcup_{k=1}^{n} S_{k})=
\sum_{k=1}^{n} vol_{d}(S_{k})
\geq n \frac{4e^{4}}{(96e^{10}\sqrt{\pi})^{d}}\  t\ vol_{d}(B_{2}^{d}) 
$$
Thus we get (2.4)
$$
vol_{d}(K \setminus K_{t}) \geq e^{-16d} n\ t\ vol_{d}(B_{2}^{d})
$$
\enddemo
\qed
\newpage

\heading
3. The Illumination Body
\endheading
\vskip 5mm

\proclaim{Theorem 3.1}
Let $K$ be a convex body in $\Bbb R^{d}$ such that
$$
\frac{1}{c_{1}}B_{2}^{d} \subset K \subset
c_{2}B_{2}^{d}
$$
Let $0 \leq t \leq (5c_{1}c_{2})^{-d-1}vol_{d}(K)$ and let $n \in \Bbb N$
with
$$
(\frac{128}{7} \pi)^{\frac{d-1}{2}} \leq n \leq 
\frac{1}{32\ edt}vol_{d}(K^{t} \setminus K)
$$
Then we have for every polytope $P_{n}$ that contains $K$ and has
at most $n$ $d-1$ dimensional faces
$$
vol_{d}(K^{t} \setminus K) \leq 10^{7}\ d^{2}(c_{1}c_{2})^{2+\frac{1}{d-1}} 
vol_{d}(P_{n} \setminus K)
$$
\endproclaim
\vskip 3mm

We want to see what this result means for bodies with a smooth
boundary. We have the  asymptotic formula [W]
$$
\lim_{t \to 0}\frac{vol_{d}(K^{t})-vol_{d}(K)}{t^{\frac{2}{d+1}}}=
\frac{1}{2}\left(\frac{d(d+1)}{vol_{d-1}(B_{2}^{d-1})}\right)^{\frac{2}{d+1}}
\int_{\partial K} \kappa(x)^{\frac{1}{d+1}}d\mu(x)
$$
Thus we get
$$
vol_{d}(K^{t})-vol_{d}(K) \sim t^{\frac{2}{d+1}}d
\int_{\partial K} \kappa(x)^{\frac{1}{d+1}} d\mu(x)
$$
And by the theorem we have
$$
n \sim \frac{1}{dt}vol_{d}(K^{t} \setminus K)
$$
Thus we get
$$
vol_{d}(K^{t})-vol_{d}(K) \sim d (\frac{1}{dn}vol_{d}(K^{t} 
\setminus K))^{\frac{2}{d+1}}\int_{\partial K} \kappa(x)^{\frac{1}{d+1}} d\mu(x)
$$
Or
$$
vol_{d}(K^{t} \setminus K)^{\frac{d-1}{d+1}} \sim d (\frac{1}{dn})^{\frac{2}{d+1}}
\int_{\partial K} \kappa(x)^{\frac{1}{d+1}} d\mu(x)
$$
$$
vol_{d}(K^{t} \setminus K) \sim d (\frac{1}{n})^{\frac{2}{d-1}}
\left(\int_{\partial K} \kappa(x)^{\frac{1}{d+1}} d\mu(x)\right)^{\frac{d+1}{d-1}}
$$
By Theorem 3.1 we get now
$$
vol_{d}(P_{n} \setminus K) \gtrsim \frac{1}{d}(\frac{1}{c_{1}c_{2}}
)^{1+\frac{d}{d+1}}(\frac{1}{n})^{\frac{2}{d-1}}
\left(\int_{\partial K} \kappa(x)^{\frac{1}{d+1}} d\mu(x)\right)^{\frac{d+1}{d-1}}
$$
By a theorem of F. John [J] we have $c_{1}c_{2} \leq d$.
\vskip 3mm

The following lemma is due to Bronshtein and Ivanov [BI] and
Dudley $[D_{1}, D_{2}]$. It can also be found in [GRS]. 
\vskip 3mm

\proclaim{Lemma 3.2}
For all dimensions $d$, $d \geq 2$, and all natural numbers $n$,
$n \geq 2d$, there is a polytope $Q_{n}$ that has $n$ vertices
and is contained in the Euclidean ball $B_{2}^{d}$ such that
$$
d_{H}(Q_{n},B_{2}^{d}) \leq \frac{16}{7}
(\frac{vol_{d-1}(\partial B_{2}^{d})}{vol_{d-1}( B_{2}^{d-1})})^
{\frac{2}{d-1}}
n^{-\frac{2}{d-1}}
$$
\endproclaim
\vskip 5mm

We have that 
$$
\align
vol_{d-1}(\partial B_{2}^{d}) & =d\ vol_{d}( B_{2}^{d})
=d \frac{\pi^{\frac{d}{2}}}{\Gamma(\frac{d}{2}+1)}        \\ 
& =d \sqrt{\pi}\frac{\Gamma(\frac{d-1}{2}+1)}{\Gamma(\frac{d}{2}+1)}
vol_{d-1}( B_{2}^{d-1}) 
\leq d \sqrt{\pi}\ vol_{d-1}( B_{2}^{d-1})
\tag 3.1
\endalign
$$
Since $d^{\frac{2}{d-1}} \leq 4$ and $(1-t)^{d} \geq 1-dt$ we get
from (3.1)
$$
d_{H}(B_{2}^{d},Q_{n}) \leq \frac{16}{7}\left(
\frac{d\sqrt{\pi}}{n}\right)^{\frac{2}{d-1}}\leq
\frac{64}{7}\pi\ n^{-\frac{2}{d-1}}
\tag 3.2
$$
\vskip 3mm

\demo{Proof of Theorem 3.1}
We denote the $d-1$-dimensional faces of $P_{n}$ by $F_{i}$, $i=1,\dots,n$,
and the cones generated by the origin and a face $F_{i}$ by $C_{i}$,
$i=1,\dots,n$. Let $x_{i} \in F_{i}$ and $\xi_{i}$, $\|\xi_{i}\|_{2}=1$,
orthogonal to $F_{i}$ and pointing to the outside of $P_{n}$. Then 
$H(x_{i},\xi_{i})$ is the hyperplane containing $F_{i}$ and 
$H^{+}(x_{i},\xi_{i})$ the halfspace containing $P_{n}$. See
figure 3.1.
We may assume that
the hyperplanes $H(x_{i},\xi_{i})$, $i=1,\dots,n$, are supporting
hyperplanes of $K$. Otherwise we can choose a polytope of
lesser volume. Let $\Delta_{i}$ be the height of the set
$$
K^{t} \cap H^{-}(x_{i},\xi_{i}) \cap C_{i}
$$
i.e. the smallest number $s$ such that
$$
K^{t} \cap H^{-}(x_{i},\xi_{i}) \cap C_{i} \subset
H^{+}(x_{i}+s \xi_{i},\xi_{i})
$$
Let $z_{i}$ be a point in $\partial K^{t} \cap C_{i}$ 
where the height $\Delta_{i}$ is
attained. 
We may assume that $B_{2}^{d} \subset K \subset cB_{2}^{d}$
where $c=c_{1}c_{2}$. Also we may assume that 
$$
P_{n} \subset 2cB_{2}^{d}      \tag 3.3
$$
if we allow twice as many faces. This follows from (3.2):
There is a polytope $Q_{k}$ such that $\frac{1}{2}B_{2}^{d} \subset
Q_{k} \subset B_{2}^{d}$ and the number of vertices $k$ is smaller than
$(\frac{128}{7}\pi)^{\frac{d-1}{2}}$. Thus $Q_{k}^{*}$ satisfies
$B_{2}^{d} \subset Q_{k}^{*} \subset 2B_{2}^{d}$ and has at most
$(\frac{128}{7}\pi)^{\frac{d-1}{2}}$ $d-1$-dimensional faces. As the new
polytope $P_{n}$ we choose the intersection of $cQ_{k}^{*}$ with the
original polytope $P_{n}$. Since we have by assumption that n is greater
than $(\frac{128}{7}\pi)^{\frac{d-1}{2}}$ the new polytope has at most
$$
\frac{1}{16\ edt}vol_{d}(K^{t} \setminus K)
\tag 3.4
$$
$d-1$-dimensional faces.

\vskip 7mm
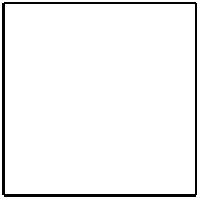
\BoxedEPSF{convex5 scaled 600}
\vskip 7mm

We show first that for $t$ with
$0 \leq t \leq (5cd)^{-d-1}vol_{d}(K)$ and all $i$, $i=1,\dots,n$
we have
$$
\Delta_{i} \leq \frac{1}{d}      \tag 3.5 
$$
Assume that there is a face $F_{i}$ with $\Delta_{i} > \frac{1}{d}$.
Consider the smallest infinite cone $D_{i}$ having $z_{i}$
as vertex and containing $K$. Since $H(x_{i},\xi_{i})$ is a supporting
hyperplane to $K$ and $K \subset c\ B_{2}^{d}$ we have
$$
K \subset D_{i} \cap H^{+}(x_{i},\xi_{i}) \cap
H^{-}(x_{i}-4c\xi_{i},\xi_{i})
$$
and
$$
D_{i} \cap H^{-}(x_{i},\xi)=[z_{i},K] \cap H^{-}(x_{i},\xi) 
$$
We have
$$
t=vol_{d}([z_{i},K] \setminus K) \geq
vol_{d}([z_{i},K] \cap H^{-}(x_{i},\xi_{i}))=
vol_{d}(D_{i} \cap H^{-}(x_{i},\xi_{i}))=
$$
$$
\frac{1}{d} \Delta_{i} vol_{d-1}(D_{i} \cap H(x_{i},\xi_{i})) \geq
\frac{1}{d^{2}}vol_{d-1}(D_{i} \cap H(x_{i},\xi_{i}))
$$
Thus
$$
vol_{d-1}(D_{i} \cap H(x_{i},\xi_{i})) \leq d^{2}t    \tag 3.6
$$
Since (3.5) does not hold we have
$$
\align
vol_{d-1}(D_{i} \cap H(x_{i}-4c\xi_{i},\xi_{i})) & =
(\frac{4c+\Delta_{i}}{\Delta_{i}})^{d-1}
vol_{d-1}(D_{i} \cap H(x_{i},\xi_{i}))                  \\ 
& \leq (4cd+1)^{d-1}vol_{d-1}((D_{i} \cap H(x_{i},\xi_{i}))
\endalign
$$
By (3.6) we get
$$
vol_{d-1}(D_{i} \cap H(x_{i}-4c\xi_{i},\xi_{i})) \leq
(4cd+1)^{d-1}d^{2}t \leq (5cd)^{d-1} d^{2} t
$$
Thus we get
$$
\align
vol_{d}(K) & \leq vol_{d}(D_{i}  
\cap  H^{+}(x_{i},\xi_{i}) 
\cap H^{-}(x_{i}-4c\xi_{i},\xi_{i}))     \\
& \leq 2c(5cd)^{d-1} d^{2} t \leq (5cd)^{d+1} t
\endalign 
$$
Thus
$$
t \geq (5cd)^{-d-1}vol_{d}(K)
$$
This is a contradiction to the assumption on $t$ in the hypothesis
of the theorem. Thus we have shown (3.5).
We consider now two cases: All those heights $\Delta_{i}$ that are
smaller than $\frac{2dt}{vol_{d-1}(F_{i})}$ and those that are greater.
We may assume that $\Delta_{i}$, $i=1,\dots,k$ are smaller than
$\frac{2dt}{vol_{d-1}(F_{i})}$ and $\Delta_{i}$, $i=k+1,\dots,n$ are
strictly greater.  
We have 
$$
vol_{d}((K^{t} \setminus P_{n}) \cap C_{i})=
\int_{0}^{\Delta_{i}} vol_{d-1}((K^{t} \setminus P_{n}) \cap C_{i}
\cap H(x_{i}+s\xi_{i},\xi_{i}))ds 
$$
Since $B_{2}^{d} \subset K \subset P_{n}$ we get
$$
vol_{d}((K^{t} \setminus P_{n}) \cap C_{i}) \leq
\int_{0}^{\Delta_{i}}vol_{d-1}(F_{i})(1+s)^{d-1}ds \leq
\Delta_{i}(1+\Delta_{i})^{d-1}vol_{d-1}(F_{i}) 
$$
By (3.5) we get
$$
vol_{d}((K^{t} \setminus P_{n}) \cap C_{i}) \leq
\Delta_{i}(1+\frac{1}{d})^{d-1}vol_{d-1}(F_{i})
$$
For $i=1,\dots,k$ we get
$$
vol_{d}((K^{t} \setminus P_{n}) \cap C_{i}) \leq
\frac{2dt}{vol_{d-1}(F_{i})} (1+\frac{1}{d})^{d-1}vol_{d-1}(F_{i}) \leq
2edt
$$
Thus we get
$$
vol_{d}((K^{t} \setminus P_{n}) \cap (\bigcup_{i=1}^{k}C_{i})) \leq
2kedt \leq 2nedt
$$
By (3.4) we get
$$
vol_{d}((K^{t} \setminus P_{n}) \cap (\bigcup_{i=1}^{k}C_{i})) \leq
\frac{1}{8}vol_{d}(K^{t} \setminus K)
\tag 3.7
$$
Now we consider the other faces. We have for $i=k+1, \dots,n$
$$
\Delta_{i} \geq \frac{2dt}{vol_{d-1}(F_{i})}
\tag 3.8
$$
We show that we have for $i=k+1, \dots,n$
$$
\Delta_{i} \leq 5c
\left(\frac{5c\ vol_{d-1}(F_{i})}{2d\ vol_{d}(K)}\right)^{\frac{1}{d-1}}
\tag 3.9
$$
Suppose that there is a face $F_{i}$ so that (3.9) does not hold. Then
we have
$$
t=vol_{d}([z_{i},K] \setminus K) \geq vol_{d}([z_{i},K] \cap
H^{-}(x_{i},\xi_{i})) = \frac{\Delta_{i}}{d}
vol_{d-1}([z_{i},K] \cap H(x_{i},\xi_{i}))
$$
Therefore we get by (3.8)
$$
vol_{d-1}([z_{i},K] \cap H(x_{i},\xi_{i})) \leq
\frac{dt}{\Delta_{i}} \leq \frac{1}{2}vol_{d-1}(F_{i})
\tag 3.10
$$
By (3.3) we have that
$$
K \subset D_{i} \cap H^{+}(x_{i},\xi_{i}) \cap H^{-}(x_{i}-4c\xi_{i},
\xi_{i})
$$
Thus
$$
vol_{d}(K) \leq vol_{d}(D_{i}  \cap H^{-}(x_{i}-4c\xi_{i},\xi_{i})) 
$$
The cone $D_{i}  \cap H^{-}(x_{i}-4c\xi_{i},\xi_{i})$ has a height equal to
$4c+\Delta_{i}$. Therefore we get
$$
vol_{d}(K) \leq
\frac{1}{d}(4c+\Delta_{i})(\frac{4c+\Delta_{i}}{\Delta_{i}})^{d-1}
vol_{d-1}(D_{i}  \cap H(x_{i},\xi_{i}))
$$
By (3.5) we have $\Delta_{i} \leq 1$. Therefore we get
$$
\align
vol_{d}(K) & \leq \frac{5c}{d}(\frac{5c}{\Delta_{i}})^{d-1}
vol_{d-1}(D_{i} \cap H(x_{i},\xi_{i}))            \\
& =\frac{5c}{d}(\frac{5c}{\Delta_{i}})^{d-1}
vol_{d-1}([z_{i},K] \cap H(x_{i},\xi_{i}))
\endalign
$$
By (3.10) we get
$$
vol_{d}(K) \leq \frac{5c}{2d}(\frac{5c}{\Delta_{i}})^{d-1}
vol_{d-1}(F_{i})
$$
This inequality implies (3.9).
\par 
Let $y_{i}$ be the unique point
$$
y_{i}=[0,z_{i}] \cap H(x_{i},\xi_{i})
$$
We want to make sure that
$y_{i} \in F_{i} \cap [z_{i},K]$. This holds since $z_{i} \in C_{i} 
\cap H^{-}(x_{i},\xi_{i})$
and $\Delta_{i}>0$. 
Since $y_{i} \in F_{i}$ we have
$$
vol_{d-1}(F_{i})=\frac{vol_{d-1}(B_{2}^{d-1})}{vol_{d-2}
(\partial B_{2}^{d-1})} \int_{\partial B_{2}^{d-1}}
r_{i}(\eta)^{d-1}d\mu(\eta)
$$
where $r_{i}(\eta)$ is the distance of $y_{i}$ to the boundary $\partial F_{i}$
in direction $\eta$, $\eta \in \partial B_{2}^{d-1}$, and, since
$y_{i} \in F_{i} \cap [z_{i},K]$, we have
$$
vol_{d-1}(F_{i} \cap [z_{i},K])=\frac{vol_{d-1}(B_{2}^{d-1})}{vol_{d-2}
(\partial B_{2}^{d-1})} \int_{\partial B_{2}^{d-1}}
\rho_{i}(\eta)^{d-1}d\mu(\eta)
$$
where $\rho_{i}(\eta)$ is the distance of $y_{i}$ to the boundary
$\partial (F_{i} \cap [z_{i},K])$. 
Consider the set
$$
A_{i}=\{ \eta |\ (1-\frac{1}{4d})r_{i}(\eta) \leq \rho_{i}(\eta)\ \}
$$
We show that
$$
\frac{1}{4} vol_{d-1}(F_{i}) \leq
\frac{vol_{d-1}(B_{2}^{d-1})}{vol_{d-2}
(\partial B_{2}^{d-1})} \int_{A_{i}^{c}}
r_{i}(\eta)^{d-1}-\rho_{i}(\eta)^{d-1}d\mu(\eta) 
\tag 3.11
$$
We have
$$
\align
\frac{vol_{d-1}(B_{2}^{d-1})}{vol_{d-2}
(\partial B_{2}^{d-1})} &  \int_{A_{i}}
r_{i}(\eta)^{d-1}-\rho_{i}(\eta)^{d-1}d\mu(\eta)    \\ 
& \leq
\frac{vol_{d-1}(B_{2}^{d-1})}{vol_{d-2}
(\partial B_{2}^{d-1})} \int_{A_{i}}
r_{i}(\eta)^{d-1}(1-(1-\frac{1}{4d})^{d-1})d\mu(\eta)  \\ 
& \leq
\frac{1}{4}\frac{vol_{d-1}(B_{2}^{d-1})}{vol_{d-2}
(\partial B_{2}^{d-1})} \int_{A_{i}}
r_{i}(\eta)^{d-1}d\mu(\eta) \leq
\frac{1}{4}vol_{d-1}(F_{i})
\endalign
$$
Therefore we get
$$
\align
\frac{vol_{d-1}(B_{2}^{d-1})}{vol_{d-2}
(\partial B_{2}^{d-1})} & \int_{A_{i}^{c}}
r_{i}(\eta)^{d-1}-\rho_{i}(\eta)^{d-1}d\mu(\eta) \geq     \\
& \frac{vol_{d-1}(B_{2}^{d-1})}{vol_{d-2}
(\partial B_{2}^{d-1})} \int_{\partial B_{2}^{d-1}}
r_{i}(\eta)^{d-1}-\rho_{i}(\eta)^{d-1}d\mu(\eta)-             \\
& \frac{vol_{d-1}(B_{2}^{d-1})}{vol_{d-2}
(\partial B_{2}^{d-1})}\int_{A_{i}}r_{i}(\eta)^{d-1}-\rho_{i}(\eta)^{d-1}d\mu(\eta) \geq  \\
& vol_{d-1}(F_{i})-vol_{d-1}(F_{i} \cap [z_{i},K])-\frac{1}{4} vol_{d-1}(F_{i})
\endalign
$$
By (3.10) we get that this is greater than $\frac{1}{4} vol_{d-1}(F_{i})$.
This implies
$$
\frac{1}{4} vol_{d-1}(F_{i}) \leq
\frac{vol_{d-1}(B_{2}^{d-1})}{vol_{d-2}
(\partial B_{2}^{d-1})} \int_{A_{i}^{c}}
r_{i}(\eta)^{d-1}-\rho_{i}(\eta)^{d-1}d\mu(\eta) 
$$
Thus we have established (3.11).
\par
We shall show that 
$$
vol_{d}((K^{t} \setminus P_{n}) \cap C_{i}) \leq 20480\ ed^{2}c^{2+\frac{1}{d-1}}
vol_{d}((P_{n} \setminus K) \cap C_{i})
\tag 3.12 
$$
We have
$$
vol_{d}(D_{i}^{c} \cap H^{+}(x_{i},\xi_{i}) \cap C_{i}) \leq 
vol_{d}((P_{n} \setminus K) \cap C_{i})
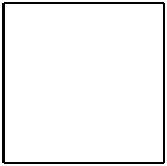
$$
\vskip 5mm

\BoxedEPSF{convex6 scaled 500}
\vskip 5mm

Compare figure 3.2. Therefore, if we want to verify (3.12) it is enough
to show
$$
vol_{d}((K^{t} \setminus P_{n}) \cap C_{i}) \leq 20480\ ed^{2}c^{2+\frac{1}{d-1}}
vol_{d}(D_{i}^{c} \cap H^{+}(x_{i},\xi_{i}) \cap C_{i}) 
$$
We may assume that $y_{i}$ and $z_{i}$ are orthogonal to $H(x_{i},\xi_{i})$.
This is accomplished by a linear, volume preserving map: Any vector
orthogonal to $\xi_{i}$ is mapped onto itself and $y_{i}$ is mapped to
$<\xi_{i},y_{i}>\xi_{i}$. See figure 3.3.
\vskip 5mm

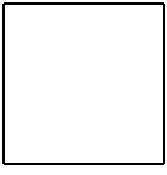
\BoxedEPSF{convex7 scaled 500}
\vskip 5mm

Let $w_{i}(\eta) \in D_{i}^{c} \cap H^{+}(x_{i},\xi_{i}) \cap C_{i}$ such that 
$w_{i}(\eta)$ is an element of the $2$-dimensional subspace containing
$0$, $y_{i}$, and $y_{i}+\eta$. Let $\delta_{i}(\eta)$ be the distance of
$w_{i}(\eta)$ to the plane $H(x_{i},\xi_{i})$. Then we have
$$
\frac{1}{d}\frac{vol_{d-1}(B_{2}^{d-1})}{vol_{d-2}
(\partial B_{2}^{d-1})} \int_{A_{i}^{c}}
(r_{i}(\eta)^{d-1}-\rho_{i}(\eta)^{d-1}) \delta_{i}(\eta)\ d\mu(\eta) \leq
vol_{d}(D_{i}^{c} \cap H^{+}(x_{i},\xi_{i}) \cap C_{i})
$$
Thus, in order to verify (3.12), it suffices to show
$$
\align
vol_{d} & ((K^{t} \setminus P_{n}) \cap C_{i}) \leq           \\ 
& 20480\ ed^{2}c^{2+\frac{1}{d-1}}
\frac{1}{d}\frac{vol_{d-1}(B_{2}^{d-1})}{vol_{d-2}
(\partial B_{2}^{d-1})} \int_{A_{i}^{c}}
(r_{i}(\eta)^{d-1}-\rho_{i}(\eta)^{d-1})\delta_{i}(\eta)\ d\mu(\eta)
\tag{3.13}
\endalign 
$$
In order to do this  
we shall show that for all $i=k+1,\dots,n$ and all $\eta \in A_{i}^{c}$ 
there is $w_{i}(\eta)$ such that the distance
$\delta_{i}(\eta)$ of $w_{i}$ from $H(x_{i},\xi_{i})$ satisfies
$$
\frac{\Delta_{i}}{\delta_{i}} \leq \left\{
\aligned
& 32dc  \hskip 50mm  \text{if} \quad 0 \leq \alpha_{i} \leq \frac{\pi}{4}     \\
& \frac{160\ dc^{2}}{r_{i}} \left(\frac{5c\ vol_{d-1}(F_{i})}
{2d\ vol_{d}(K)} \right)^{\frac{1}{d-1}}  \hskip 12mm
\text{if} \quad \frac{\pi}{4} \leq \alpha_{i} \leq \frac{\pi}{2}  
\endaligned
\right .
\hskip 7mm
\tag 3.14
$$
The angles $\alpha_{i}(\eta)$ and $\beta_{i}(\eta)$ are given in figure 3.3.
We have for all $\eta \in A_{i}^{c}$
$$
\aligned
\delta_{i}=& (r_{i}-\rho_{i})\frac{\sin(\alpha_{i})\sin(\beta_{i})}
{\sin(\pi-\alpha_{i}-\beta_{i})}                                         \\
\Delta_{i}=& \rho_{i} \tan \alpha_{i}
\endaligned
\hskip 10mm
0 \leq \alpha_{i},\beta_{i} \leq \frac{\pi}{2}
\tag 3.15
$$
Thus we get
$$
\frac{\Delta_{i}}{\delta_{i}} \leq
\frac{\rho_{i}}{r_{i}-\rho_{i}}
\frac{\sin(\pi-\alpha_{i}-\beta_{i})}{\cos(\alpha_{i}) \sin(\beta_{i})}
\leq \frac{\rho_{i}}{(r_{i}-\rho_{i})\cos(\alpha_{i}) \sin(\beta_{i})}
$$
By (3.11) we have $\rho_{i} \leq (1-\frac{1}{4d})r_{i}$. Therefore
we get
$$
\frac{\Delta_{i}}{\delta_{i}} \leq
\frac{4d}{\cos(\alpha_{i}) \sin(\beta_{i})}
$$
Since $B_{2}^{d} \subset K \subset P_{n} \subset 2c\ B_{2}^{d}$
we get that $\tan \beta_{i} \geq \frac{1}{4c}$: Here we have to take into account that
we applied a transform to $K$ mapping $y_{i}$ to $<\xi_{i},y_{i}>\xi_{i}$.
That leaves the distance of $F_{i}$ to the origin unchanged and $r_{i}(\eta)$
is less than $4c$. 
If $\beta_{i} \geq \frac{\pi}{4}$
we have $\sin \beta_{i} \geq \frac{1}{\sqrt 2}$. If $\beta_{i} \leq \frac{\pi}{4}$
then $\frac{1}{4c} \leq \tan \beta_{i} =\frac{\sin \beta_{i}}{\cos \beta_{i}} \leq
\sqrt{2} \sin \beta_{i}$. Therefore we get
$$
\frac{\Delta_{i}}{\delta_{i}} \leq \frac{16 \sqrt{2}\ dc}{\cos \alpha}
$$
Therefore we get for all $0 \leq \alpha_{i} \leq \frac{\pi}{4}$
$$
\frac{\Delta_{i}}{\delta_{i}} \leq 32dc
$$
By (3.9) and (3.15) we get
$$
\frac{\Delta_{i}}{\delta_{i}} \leq \frac{1}{r_{i}-\rho_{i}}
\frac{\sin(\pi-\alpha_{i}-\beta_{i})}{\sin(\alpha_{i}) \sin(\beta_{i})}
5c\left(\frac{5c\ vol_{d-1}(F_{i})}
{2d\ vol_{d}(K)} \right)^{\frac{1}{d-1}}
$$
We proceed as in the estimate above and obtain
$$
\frac{\Delta_{i}}{\delta_{i}} \leq \frac{16 \sqrt{2} \ dc}{r_{i}}
\frac{5c}{\sin(\alpha_{i})}
\left(\frac{5c\ vol_{d-1}(F_{i})}
{2d\ vol_{d}(K)} \right)^{\frac{1}{d-1}}
$$
Thus we get for $\frac{\pi}{4} \leq \alpha_{i} \leq \frac{\pi}{2}$
$$
\frac{\Delta_{i}}{\delta_{i}} \leq \frac{32 \ dc}{r_{i}}
5c \left(\frac{5c\ vol_{d-1}(F_{i})}
{2d\ vol_{d}(K)} \right)^{\frac{1}{d-1}}
$$
We verify now (3.13). By the definition of $A_{i}$ we get
$$
\align
\frac{vol_{d-1}(B_{2}^{d-1})}{vol_{d-2}
(\partial B_{2}^{d-1})}  & \int_{A_{i}^{c}}  
(r_{i}(\eta)^{d-1}-\rho_{i}(\eta)^{d-1})\delta_{i}(\eta)\ d\mu(\eta) \geq   \\
& (1-e^{-\frac{1}{8}})\frac{vol_{d-1}(B_{2}^{d-1})}{vol_{d-2}
(\partial B_{2}^{d-1})} \int_{A_{i}^{c}}r_{i}(\eta)^{d-1}\delta_{i}\ d\mu(\eta)
\endalign
$$
We get by (3.15)
$$                           
\frac{1}{320dc} \Delta_{i} \frac{vol_{d-1}(B_{2}^{d-1})}
{vol_{d-2}(\partial B_{2}^{d-1})}
\left\{\int_{A_{i}^{c} \atop \alpha_{i} \leq \frac{\pi}{4}}r_{i}^{d-1}d\mu+
\frac{1}{5c}\left(\frac{2d\ vol_{d}(K)}{5c\ vol_{d-1}(F_{i})}\right)^{\frac{1}{d-1}}
\int_{A_{i}^{c} \atop {\alpha_{i} >\frac{\pi}{4}}} r_{i}^{d}d\mu \right\}
$$
By (3.11) we get that either
$$
\frac{vol_{d-1}(B_{2}^{d-1})}{vol_{d-2}(\partial B_{2}^{d-1})}
\int_{A_{i}^{c} \atop \alpha_{i} \leq \frac{\pi}{4}}r_{i}^{d-1}d\mu
\geq \frac{1}{8}\ vol_{d-1}(F_{i})
$$
or
$$
\frac{vol_{d-1}(B_{2}^{d-1})}{vol_{d-2}(\partial B_{2}^{d-1})}
\int_{A_{i}^{c} \atop \alpha_{i} > \frac{\pi}{4}}r_{i}^{d-1}d\mu
\geq \frac{1}{8}\ vol_{d-1}(F_{i})
$$
In the first case we get for the above estimate
$$
\align
\frac{vol_{d-1}(B_{2}^{d-1})}{vol_{d-2}
(\partial B_{2}^{d-1})}  \int_{A_{i}^{c}}  &  
(r_{i}(\eta)^{d-1}-\rho_{i}(\eta)^{d-1})\delta_{i}(\eta)\ d\mu(\eta) \geq   \\
& \frac{\Delta_{i}}{2560dc}\ vol_{d-1}(F_{i}) \geq \frac{1}{2560edc}\ 
vol_{d}((K^{t} \setminus P_{n}) \cap C_{i})
\endalign
$$
The last inequality is obtained by using (3.5): Since $B_{2}^{d} \subset
K$ we have for all hyperplanes $H$ that are parallel to $F_{i}$ 
$vol_{d-1}(K^{t} \cap H \cap C_{i}) \leq (1+\Delta_{i})^{d-1} vol_{d-1}(F_{i})$.
By (3.5) we get $vol_{d-1}(K^{t} \cap H \cap C_{i}) \leq e\ vol_{d-1}(F_{i})$.
In the second case we have
$$
\align
\frac{vol_{d-1}(B_{2}^{d-1})}{vol_{d-2}
(\partial B_{2}^{d-1})}  & \int_{A_{i}^{c}}  
(r_{i}(\eta)^{d-1}-\rho_{i}(\eta)^{d-1})\delta_{i}(\eta)\ d\mu(\eta) \geq      \\
& \frac{1}{5c}\left(\frac{2d\ vol_{d}(K)}{5c\ vol_{d-1}(F_{i})}\right)^{\frac{1}{d-1}}
\frac{1}{320dc} \Delta_{i} \frac{vol_{d-1}(B_{2}^{d-1})}
{vol_{d-2}(\partial B_{2}^{d-1})}
\int_{A_{i}^{c} \atop \alpha_{i} > \frac{\pi}{4}}r_{i}^{d}d\mu \geq       \\
& \frac{1}{5c}\left(\frac{2d\ vol_{d}(K)}{5c\ vol_{d-1}(F_{i})}\right)^{\frac{1}{d-1}}
\frac{1}{320dc} \Delta_{i} \frac{vol_{d-1}(B_{2}^{d-1})}
{(vol_{d-2}(\partial B_{2}^{d-1}))^{\frac{d}{d-1}}}
\left(\int_{A_{i}^{c} \atop \alpha_{i} > \frac{\pi}{4}}r_{i}^{d-1}d\mu 
\right)^{\frac{d}{d-1}} \geq                                               \\
& \frac{1}{5c}\left(\frac{2d\ vol_{d}(K)}{5c\ vol_{d-1}(F_{i})}\right)^{\frac{1}{d-1}}
\frac{\Delta_{i}}{320dc} vol_{d-1}(B_{2}^{d-1})^{-\frac{1}{d-1}}
(\frac{1}{8}vol_{d-1}(F_{i}))^{\frac{d}{d-1}}=    \\
& \frac{1}{5c}\left(\frac{d\ vol_{d}(K)}{20c\ vol_{d-1}(B_{2}^{d-1})}
\right)^{\frac{1}{d-1}}\frac{\Delta_{i}}{2560dc}\ vol_{d-1}(F_{i}) \geq     \\
& \frac{1}{5c}\left(\frac{d\ vol_{d}(K)}{20c\ vol_{d-1}(B_{2}^{d-1})}
\right)^{\frac{1}{d-1}}\frac{1}{2560edc}\ vol_{d}
((K^{t} \setminus P_{n}) \cap C_{i})  
\endalign
$$
Since $B_{2}^{d} \subset K$ we get
$$
\align
\frac{vol_{d-1}(B_{2}^{d-1})}{vol_{d-2}
(\partial B_{2}^{d-1})}  & \int_{A_{i}^{c}}  
(r_{i}(\eta)^{d-1}-\rho_{i}(\eta)^{d-1})\delta_{i}(\eta)\ d\mu(\eta) \geq          \\
& \frac{1}{5c}\left(\frac{d\ vol_{d}(B_{2}^{d})}{20c\ vol_{d-1}(B_{2}^{d-1})}
\right)^{\frac{1}{d-1}}\frac{1}{2560edc}\ vol_{d}((K^{t} 
\setminus P_{n}) \cap C_{i}) \geq  \\
& \frac{1}{5c}\left(\frac{1}{20c}
\right)^{\frac{1}{d-1}}\frac{1}{2560edc}\ vol_{d}((K^{t} 
\setminus P_{n}) \cap C_{i})  \geq                             \\           
& (20480\ edc^{2+\frac{1}{d-1}})^{-1}\ vol_{d}((K^{t} 
\setminus P_{n}) \cap C_{i})  
\endalign
$$
The second case gives a weaker estimate. Therefore we get for both cases
$$
\align
vol_{d}((K^{t} & \setminus P_{n}) \cap C_{i}) \leq     \\
& 20480\ edc^{2+\frac{1}{d-1}}
\frac{vol_{d-1}(B_{2}^{d-1})}{vol_{d-2}
(\partial B_{2}^{d-1})}  \int_{A_{i}^{c}}  
(r_{i}(\eta)^{d-1}-\rho_{i}(\eta)^{d-1})\delta_{i}\ d\mu(\eta)
\endalign  
$$
Thus we have verified (3.13) and by this also (3.12). By (3.12) we get
$$
\align
vol_{d}((K^{t} \setminus P_{n}) \cap (\bigcup_{i=k+1}^{n}C_{i})) &\leq
20480\ ed^{2}c^{2+\frac{1}{d-1}}\ vol_{d}((\bigcup_{i=k+1}^{n}C_{i}) \cap (P_{n}
\setminus K))  \\
& \leq 20480\ ed^{2}c^{2+\frac{1}{d-1}}\ vol_{d}( (P_{n} \setminus K))  
\tag 3.16
\endalign                                                      
$$
If the assertion of the theorem does not hold we have 
$$
vol_{d}( (P_{n} \setminus K)) \leq \frac{1}{8} \frac{vol_{d}(K^{t} \setminus K)}
{20480\ ed^{2}c^{2+\frac{1}{d-1}}}                                     \tag 3.17
$$
Thus we get
$$
vol_{d}((K^{t} \setminus P_{n}) \cap (\bigcup_{i=k+1}^{n}C_{i})) \leq
\frac{1}{8} vol_{d}(K^{t} \setminus K)
$$
Together with (3.7) we obtain
$$
vol_{d}(K^{t} \setminus P_{n}) \leq \frac{1}{4} vol_{d}(K^{t} \setminus K)
\leq \frac{1}{4}\{vol_{d}(K^{t} \setminus P_{n})+vol_{d}(P_{n} \setminus K)\}
\tag 3.18
$$
By (3.17) we have
$$
\align
vol_{d}(P_{n} \setminus K)  & \leq \frac{1}{8} \frac{vol_{d}(K^{t} \setminus K)}
{20480\ ed^{2}c^{2+\frac{1}{d-1}}}               \\  
& \leq \frac{1}{2} vol_{d}(K^{t} \setminus K)
\leq \frac{1}{2}vol_{d}(K^{t} \setminus P_{n})+ \frac{1}{2}vol_{d}
(P_{n} \setminus K)        
\endalign
$$
This implies
$$
vol_{d}(P_{n} \setminus K) \leq vol_{d}(K^{t} \setminus P_{n})
$$
Together with (3.18) we get now the contradiction
$$
vol_{d}(K^{t} \setminus P_{n}) \leq 
\frac{1}{2} vol_{d}(K^{t} \setminus P_{n})
$$    
\enddemo
\qed

\enddocument
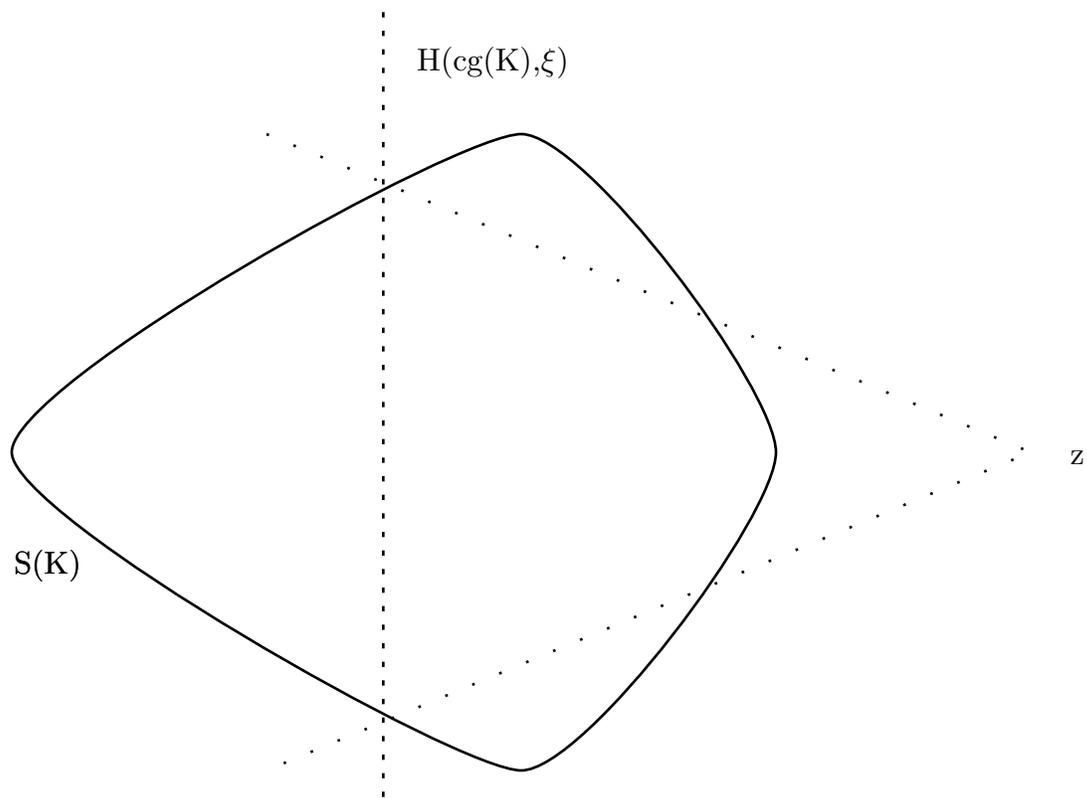
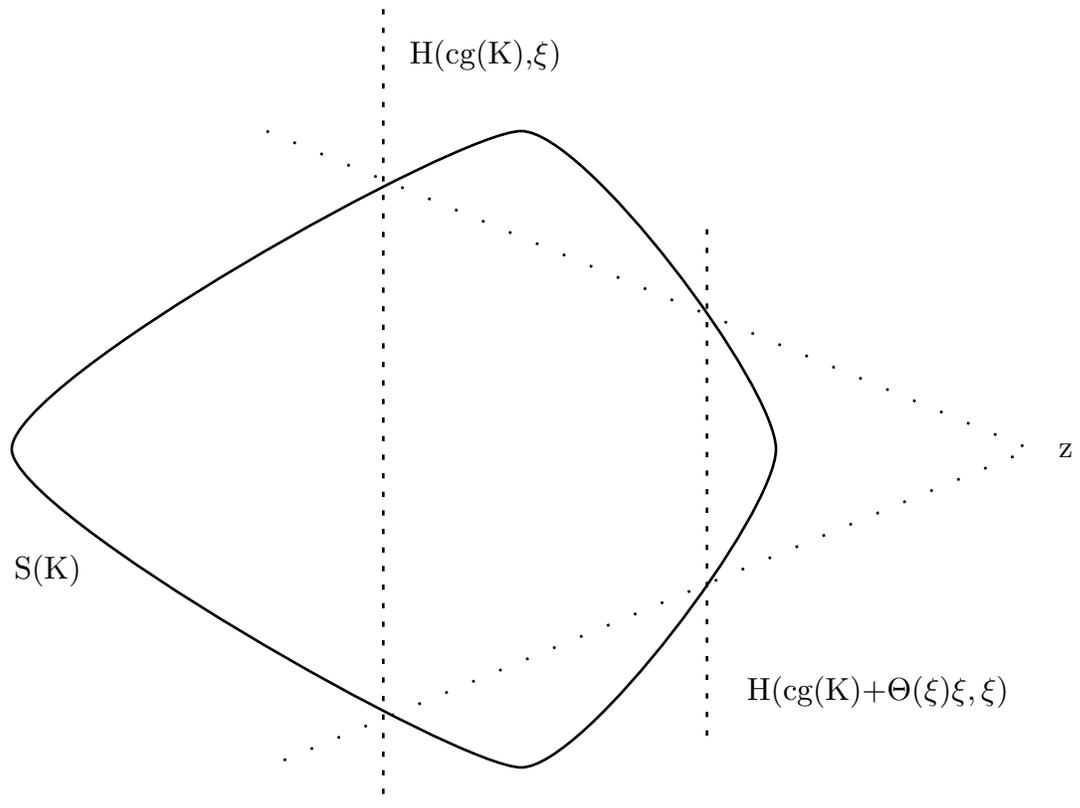
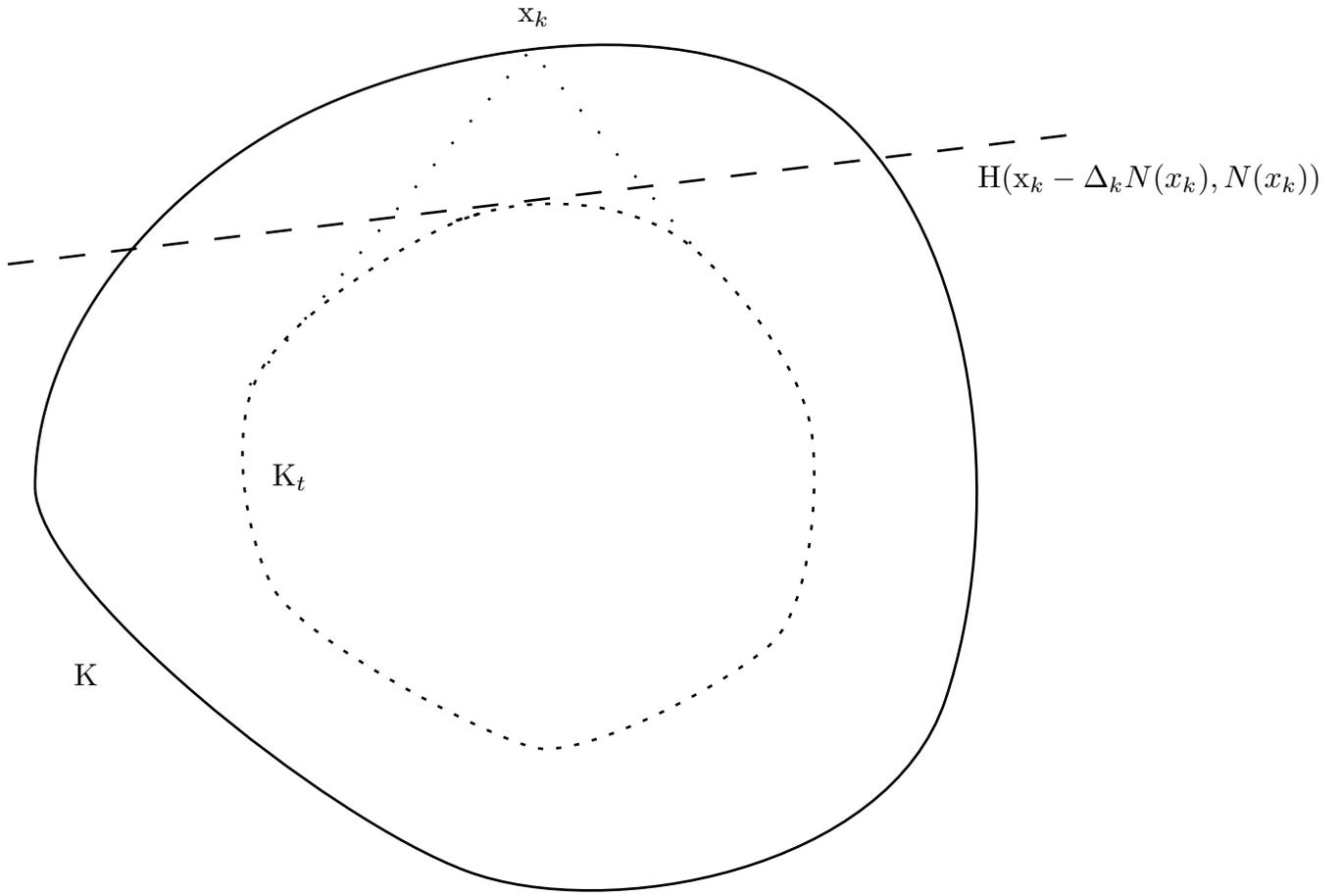
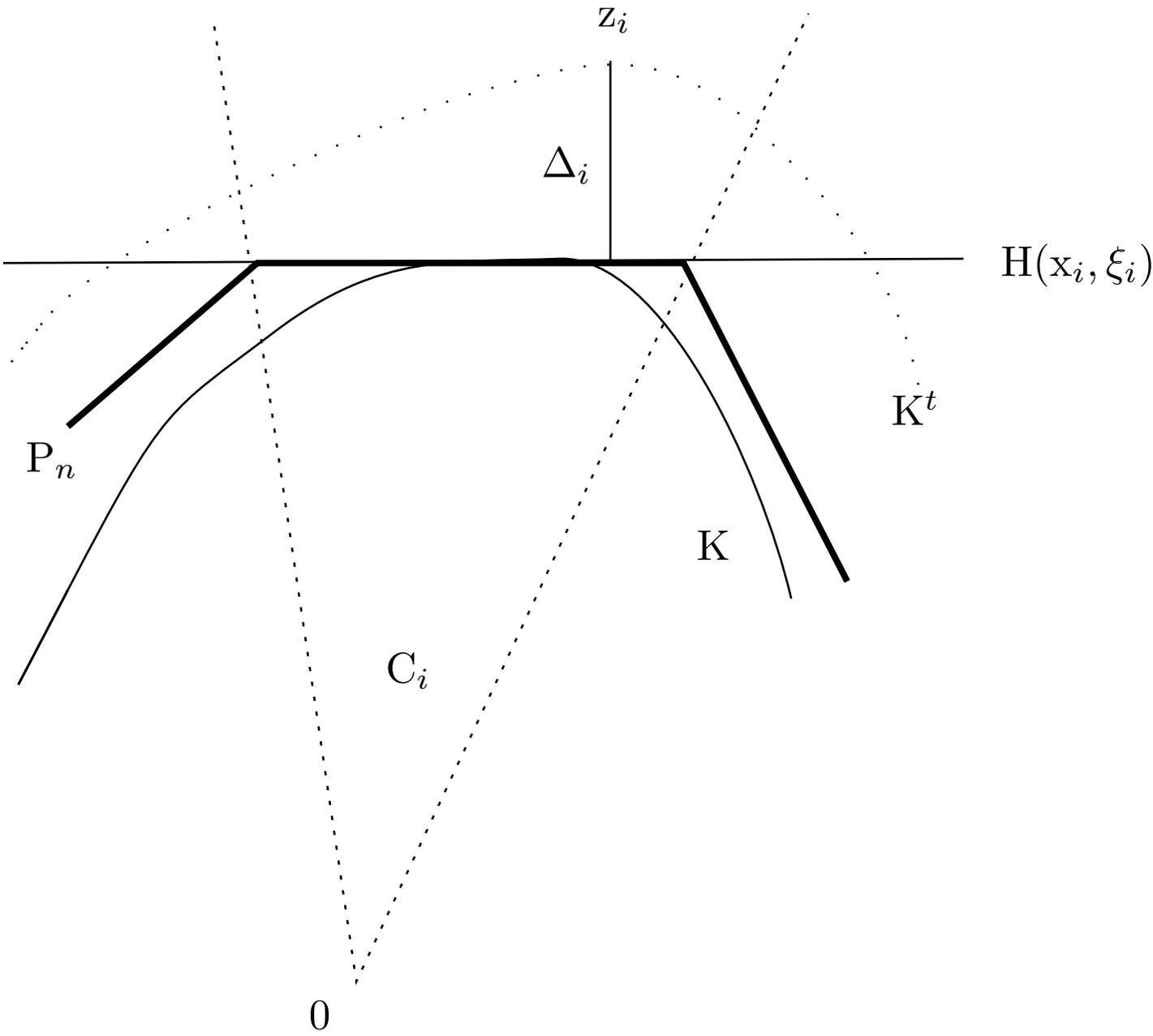
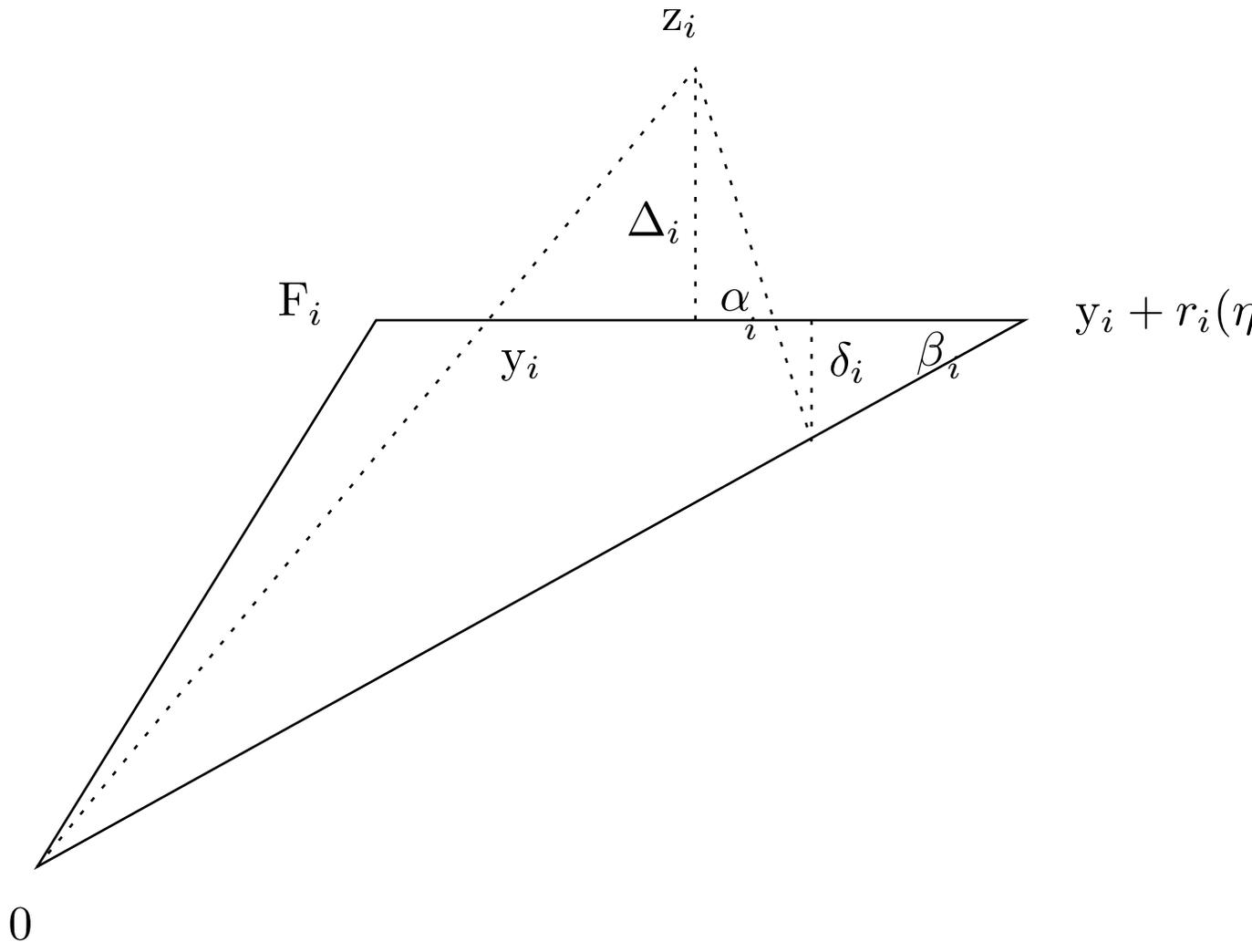
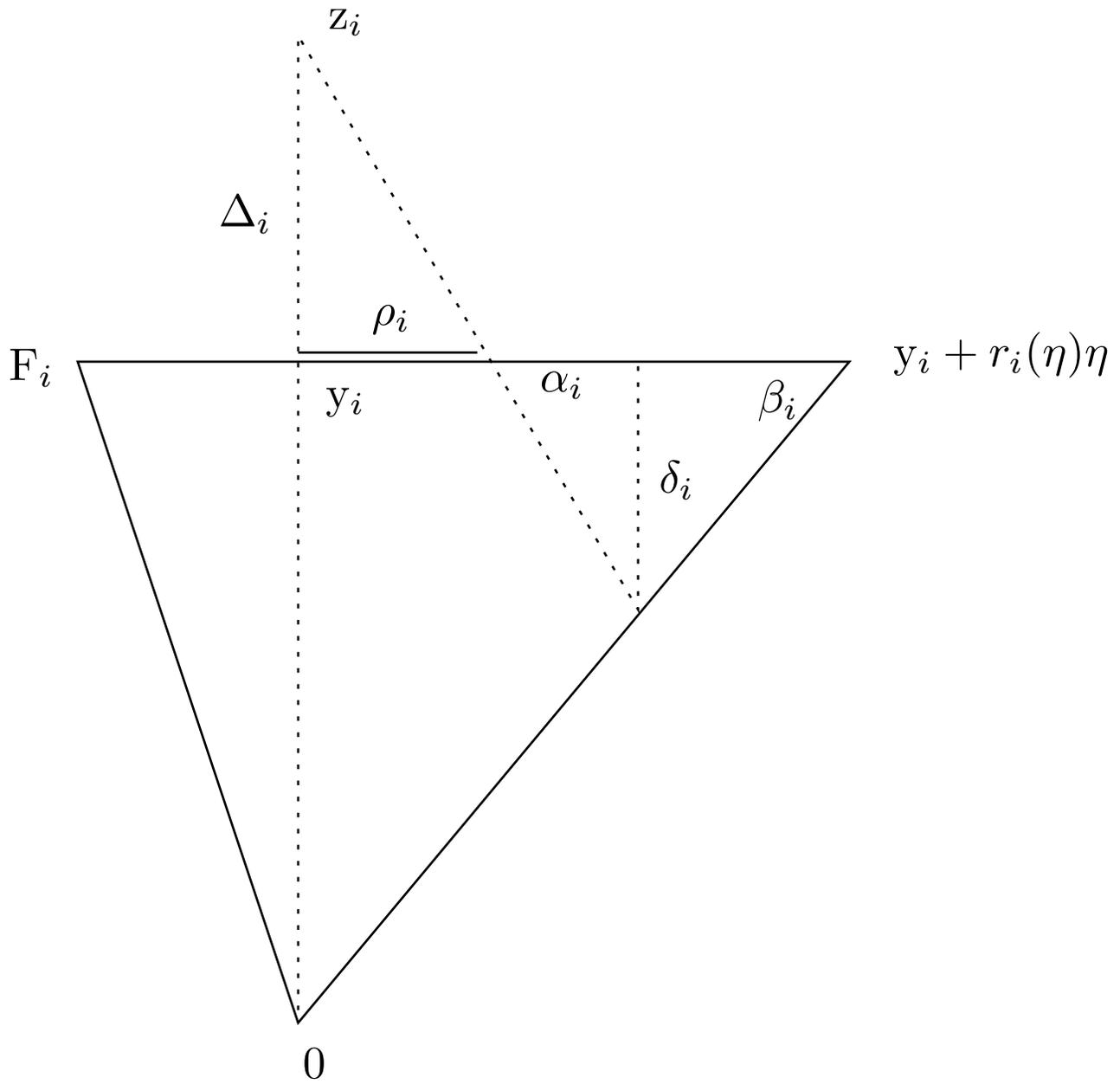
\bye